\documentclass[12pt]{article} 
\usepackage{amssymb,latexsym}
\def\neweq{\setcounter{equation}{0}}
\newtheorem{theorem}{Theorem}[section]
\newtheorem{pr}[theorem]{Proposition}
\newtheorem{cor}[theorem]{Corollary}
\newtheorem{de}[theorem]{Definition}
\newtheorem{con}[theorem]{Conjecture}
\newtheorem{rem}[theorem]{Remark}
\newtheorem{lem}[theorem]{Lemma}

\newtheorem{problem}[theorem]{Problem}
\newtheorem{ex}[theorem]{Example}

\font\germ=eufm10

\def\Q{\mathbb{Q}}

\def\A{{\cal A}}
\def\B{{\cal B}}
\def\D{{\cal D}}
\def\E{{\cal E}}
\def\G{{\cal G}}
\def\H{{\cal H}}
\def\I{{\cal I}}
\def\K{{\cal K}}
\def\L{{\cal L}}
\def\M{{\cal M}}
\def\P{{\cal P}}
\def\R{{\cal R}}

\def\s{\hbox{\germ S}}
\def\<{\langle}
\def\>{\rangle}

\def\Z{\mathbb{Z}}
\def\C{\mathbb{C}}
\def\wt{\widetilde}
\def\wh{\widehat}

\def\ds{\displaystyle}
\def\qed{\hfill$\vrule height 2.5mm width 2.5mm depth 0mm$}

\title{On some quadratic algebras
\date{ }
\author{Anatol N. Kirillov  \\~ \\
{\small {\it CRM, University of Montreal}} \\
{\small {\it C.P. 6128, Succursale A, Montreal (Quebec), H3C 3J7, 
Canada}}\\
{\small {\it and}}
 \\
{\small {\it Steklov Mathematical Institute,}} \\
{\small {\it Fontanka 27, St.Petersburg, 191011, Russia}}}
}
\begin{document}

\maketitle

\begin{abstract}
We study some quadratic 
algebras which are appeared in the low--dimen\-sio\-nal topology and Schubert 
calculus.
We introduce the Jucys--Murphy elements in the braid algebra and in the pure 
braid group, as well as the Dunkl elements in the extended affine braid group. 
Relationships between the Dunkl elements, Dunkl operators and 
Jucys--Murphy elements are described. 
\end{abstract}

\section{Introduction}
\label{sec:intro}
\neweq

In this paper we study some quadratic algebras which are naturally 
appeared in the low--dimensional topology [B-N], [Dr2], [Ko2] and Schubert 
calculus [FK], and investigate some of their properties. We define a quadratic 
algebra $\G_n$, which is a further generalization of the quadratic algebras 
$\E_n$ and $\E_n^t$ introduced in [FK], 
Sections~2 and 15. Our main idea is to apply the results obtained for the 
braid algebra $\B_n$, [Dr2], [Ko2], to the algebra $\G_n$. The main 
observation is that some results which are well--known for the braid 
algebra $\B_n$, can be 
(re)proven for the quadratic algebra $\G_n$. For example,
it happens that the quadratic algebra $\G_n$ and the braid algebra $\B_n$ 
have the same Hilbert series 
$$H(\B_n;t)=H(\G_n;t)=\frac{1}{(1-t)(1-2t)\cdots (1-(n-1)t)},
$$
%\vbox{\hskip -0.6cm 
and have the commutative subalgebras $\K_n\subset\B_n$ and 
$\H_n\subset\G_n$ which are both isomorphic to the ring of polynomials.  
The algebra $\K_n$ is generated by the 
Jucys--Murphy elements, whereas the algebra $\H_n$ is generated by the Dunkl 
ones. Both algebras $\G_n$ and $\B_n$ have the additive bases consisting 
of all monomials in normal form (Theorem~2.4 and Corollary~4.3; cf. [L]). %}
We expect also that the certain quotients $\E_n^0$ and $\B_n^0$ of the 
quadratic algebras $\G_n$ and $\B_n$ respectively, have the same Hilbert 
polynomials as well (see Section~9 for details).

One of the main goals of this paper is to describe the 
relationships between the Dunkl and Jucys--Murphy elements, and to compute 
the dimensions of 
homogeneous component of degree $\le 6$ of the quadratic algebra $\E_n^0$.

The structure of the paper is following:

In Section~2 we define the braid algebra $\B_n$ as the infinitesimal 
deformation of the pure braid group $P_n$ (cf. T.~Kohno [Ko2]).
%It is interesting to note that 
The algebra $\B_n$ and its completion $\wh\B_n$ have many important 
applications to the low--dimensional topology, see e.g. [B-N], [Dr2], [F], 
[Ko3], [L].

In Section~3 we define the Jucys--Murphy elements $d_k$ in the braid 
algebra $\B_n$, and the multiplicative Jucys--Murphy elements $D_k$ in 
the pure braid group $P_n$. Follow to A.~Ram, [Ra], we prove that the 
Jucys--Murphy element $d_k\in\B_n$ is the quasi--classical limit of the 
element $D_k\in P_n$. We prove also, that the infinitesimal deformation 
of the multiplicative Jucys--Murphy element $D_k$ coincides with the 
element $d_k$.

In Section~4 we define the 
quadratic algebra $\G_n$ and compute the Hilbert series for this algebra.
In Section~5 we define the Dunkl elements $\theta_j$ in the quadratic algebra 
$\G_n$ and describe a commutative subalgebra generated by $\theta_j$, 
$1\le j\le n$. 

In Section~6, which contains one of the main results of this paper, 
we study a relationship between the Jucys--Murphy elements 
and the Dunkl elements. We define the dual Dunkl elements $Y_1^*,\ldots 
Y_n^*$ as certain elements in the extended affine braid group $\wt B_n$, and 
prove that under a natural homomorphism $\pi :\wt B_n\to B_n$, the dual 
Dunkl element $Y_k^*$ maps to the multiplicative Jucys--Murphy element 
$D_k$. We explain also a connection between Dunkl element $Y_k\in\wt B_n$ 
and the Dunkl operator ${\bf Y}_k\in\D_{q,x}[W]$, where $\D_{q,x}[W]$ is 
the algebra of $q$--difference operators with permutations.

In Section~7 we study the relations in the algebra $\E_n^t$. In 
Section~8 we use the properties of "normal basis" introduced in Section~4, 
to study the 
quotient algebra $\E_n^0$, and compute the dimensions of homogeneous 
components $\E_{n,k}^0$,
$1\le k\le 6$, and also the Hilbert polynomial $H(\E_4^0;t)$.
The last polynomial was also computed, using computer, by J.-E.~Ross [Ro], 
who computed the Hilbert polynomial $H(\E_5^0;t)$ as well.

In Section~9 we consider a certain quotient $\B_n^0$ of the braid algebra 
$\B_n$, and make a conjecture that the quadratic algebra $\E_n^0$ and 
$\B_n^0$ have the same Hilbert polynomials. In particular, we expect that 
the algebra $\B_n^0$ (as well as $\E_n^0$) has a finite dimension.

In Section~10 we study the commutative quadratic algebra $\A_n^t$ (denote by 
$\E C_n$ in [FK], Section~4.3), which is a commutative quotient of $\E_n^t$ 
(see Definition~\ref{d7.1}), and explain some details which were omitted in 
[FK].

\medskip

\textsc{Acknowledgments.}
This work was done during my stay at the CRM, University of Montreal. I would 
like to thank all my colleagues for discussions and support. I would like 
also to thank M.~Noumi and A.N.~Varchenko for fruitful and stimulating 
discussions. My 
special thanks to Arun Ram, who explained me during my stay at the 
University of Wisconsin (1995) the importance of the Jucys--Murphy 
elements in the representation theory, and N.A.~Liskova for an 
inestimable help on different stages of this paper and constructive 
criticism.

\section{Quadratic algebra $\B_n$}
\label{sec:quadB}
\neweq

We start with consideration of the quadratic algebra $\B_n$ which is the 
infinitesimal deformation of the pure braid group $P_n$, [Ko2]. This 
quadratic algebra is well--known as the algebra of (finite order) 
Vassiliev invariants of braids, [Ko3], [F], [B-N]. The completion 
$\wh\B_n$ of $\B_n$ with respect to the grading was considered in the 
papers of V.~Drinfeld, [Dr1], [Dr2], in his study of quasi--Hopf algebras, 
and in the papers of T.~Kohno, [Ko1]--[Ko3], in his study of monodromy 
representations of braid groups. 

\begin{de}\label{d3.1}{\rm([Ko2], [Dr2])} 
Define the braid algebra $\B_n$ as the quadratic algebra 
(say over $\Z$) with generators $X_{ij}$, $1\le i<j\le n$, which satisfy 
the following relations
\begin{equation}\label{3.1}
 1)~~X_{ij}\cdot X_{jk}-X_{jk}\cdot X_{ij}=X_{ik}\cdot X_{ij}-
	X_{ij}\cdot X_{ik}=X_{jk}\cdot X_{ik}-X_{ik}\cdot X_{jk};
\end{equation}
\begin{equation}\label{3.2}
2)~~X_{ij}\cdot X_{kl}=X_{kl}\cdot X_{ij},~~{\rm if~~all}~~i,j,k,l
	~~{\rm are~~distinct}.~~~~~~~~~~~~~~~~~~~
\end{equation}
\end{de}

\begin{rem} {\rm The algebra $\B_n$ is denoted by $\P_n$ in [L], and as 
$A^n$ in [Ko1]--[Ko3].}
\end{rem}
The algebra $\B_n$ has a natural structure of a cocommutative Hopf 
algebra. Namely, a comultiplication $\Delta 
:~\B_n\to\B_n\otimes\B_n$, antipode $S:~\B_n\to\B_n$, and counit 
$\epsilon :~A\to\Z$, can be define as follows
\begin{eqnarray*}
\Delta (X_{ij})&=&1\otimes X_{ij}+X_{ij}\otimes 1,\\
S(X_{ij})&=&-X_{ij},\\
\epsilon (X_{ij})&=&0,\ \ \epsilon (1)=1.
\end{eqnarray*}

Let us explain briefly an origin of the relations 
(\ref{3.1})--(\ref{3.2}). For more details and proofs see [Dr1], [Dr2], 
[F], [Ko1]--[Ko3], [L]. Let us denote by
$$\M=\C^n\setminus\bigcup_{1\le i<j\le n}\{ z_i=z_j\} ,
$$
the configuration space of $n$ distinct ordered points in $\C$. It is 
well--known that the fundamental group $\pi_1(\M )$ of the space $\M$ 
coincides with the pure braid group $P_n$:
$$\pi_1(\M )=P_n.
$$
Let
$$w_{ij}=w_{ji}=\frac{1}{2\pi\sqrt -1}d\log (z_i-z_j)
$$
be closed 1--form on $\M$.
Then, $\{ w_{ij}~|~1\le i<j\le n\}$ represents a basis for $H^1(\M )$, 
and the relations among $w_{ij}$, $1\le i<j\le n$, are generated by the 
Arnold relations
$$w_{ij}\wedge w_{jk}+w_{jk}\wedge w_{ik}+w_{ik}\wedge w_{ij}=0
$$
for $i<j<k$ (see, e.g. [A]).

Let $\{ X_{ij}=X_{ji}~|~1\le i<j\le n\}$ be a set of non--commutative 
variables, and consider a formal connection
\begin{equation}
\Omega =\sum_{1\le i<j\le n}w_{ij}X_{ij}. \label{4.4}
\end{equation}

\begin{lem}\label{l3.1} {\rm (T.~Kohno)} The connection $\Omega$ is 
integrable if and only if the following conditions are satisfied
\begin{eqnarray*}
&& [X_{ij},X_{ik}+X_{jk}]=[X_{ij}+X_{ik}, X_{jk}]=0,  \ {\rm for} \ 
i<j<k,\\
&&[X_{ij},X_{kl}]=0, \ {\rm for~ distinct} \ i,j,k,l.
\end{eqnarray*}
\end{lem}
Let us remind that $[a,b]=ab-ba$ is the usual commutator.

It is clear that these conditions are equivalent to the relations 
(\ref{3.1})--(\ref{3.2}). A proof of Lemma~\ref{l3.1} follows from the 
observations that $d\Omega =0$, and the integrability 
condition for the connection $\Omega$ is equivalent to the following one: 
$\Omega\wedge\Omega =0$. Notice that the connection (\ref{4.4}) is the 
formal version of the so--called Knizhnik--Zamolodchikov connection in 
conformal field theory.

The relations (\ref{3.1}) and (\ref{3.2}) can be also interpreted [Ch2], 
[Dr2], [Ko2], as the self--consistency conditions
$$\frac{\partial A_j}{\partial z_i}-\frac{\partial A_i}{\partial z_j}
=[A_i,A_j], \ \ 1\le i,j\le n,
$$
of the Knizhnik--Zamolodchikov (system of) equation(s)
$$\frac{\partial\Phi}{\partial z_i}=kA_i\Phi, \ \ 1\le i\le n,
$$
where $A_i=\displaystyle\sum_{j\ne i}\frac{X_{ij}}{z_i-z_j}$, and a 
function $\Phi :=\Phi (z)$, $z=(z_1,\ldots ,z_n)$, takes the values in a 
certain algebra (say over $\C((z))$) generated by the elements $X_{ij}$.

Below we formulate the basic properties of the quadratic algebra $\B_n$.

\begin{theorem}\label{t3.1} {\rm ([Ko1], [Ko2])} 
The Hilbert series of the algebra $\B_n$ is given by
\begin{equation}\label{3.3}
H(\B_n;t)=\frac{1}{(1-t)(1-2t)\cdots (1-(n-1)t)}.
\end{equation}
\end{theorem}

\begin{cor}\label{c3.2} {\rm ([Dr1], [L])} A linear basis in the algebra $\B_n$ 
is given 
by the following set of noncommutative monomials
\begin{equation}
Z_2\cdot\ldots\cdot Z_n, \label{4.5}
\end{equation}
where
\begin{eqnarray*}
Z_2 & = & \{ 1,X_{12}\},\\
& \cdots & \\
Z_k & = & \{ 1\}\cup\{ X_{i_1k}\ldots X_{i_lk}~|~1\le i_j<k,~1\le j\le 
l\}, \ \ 2\le k\le n.
\end{eqnarray*}
\end{cor}
A monomial in $X_{ij}$'s is called to be in {\it normal form} if it 
belongs to the set (\ref{4.5}).

\begin{rem} {\rm Corollary \ref{c3.2} was stated without a proof in [Dr1]. 
The formula (\ref{3.3}) was obtained in [Ko1] by constructing a free 
resolution of $\C$ as a trivial $\B_n$--module. Also, the formula 
(\ref{3.3}) is a consequence of the fact that $\B_n$ is a semi-direct 
product of free associative algebras as shown in [Dr2]; see also [L]. 
Note, that the relations (\ref{3.1}) can be interpreted as the horizontal 
version of the 4--term relations in the theory of Vassiliev link 
invariants.}
\end{rem}

Let us consider a completion $\widehat\B_n$ of the algebra 
$\B_n\otimes_{\Z}\C$ with respect to the powers of the ideal $\I=(X_{ij})$ 
generated by $X_{ij}$, $1\le i<j\le n$. More precisely, let us denote by 
$\C\langle X_{ij}\rangle$ the ring of non--commutative formal power 
series with indeterminates $X_{ij}$, $1\le i<j\le n$, and let $J$ be an 
ideal generated by the relations (\ref{3.1}) and (\ref{3.2}). Then $\widehat 
\B_n=\C\langle X_{ij}\rangle /J$. It is well--known [Ko1], that 
there exists an isomorphism of complete Hopf algebras
$$\widehat{\C [P_n]}\cong \widehat\B_n,
$$
where $\widehat{\C [P_n]}$ stands for the completion of the group ring 
of the pure braid group $P_n$ with respect to the powers of the 
augmentation ideal. %For further details see [Ko1], [Ko2], [Dr2], [B-N].

\section{Jucys--Murphy elements}
\label{sec:Juc}
\neweq

\begin{de}\label{d4.1} The Jucys--Murphy elements $d_j$, for $j=2,\ldots 
,n$, in the quadratic algebra $\B_n$ are defined by
\begin{equation}\label{4.1}
d_j=\sum_{1\le i<j}X_{ij}.
\end{equation}
\end{de}
It is clear that $d_i$ is a primitive element, i.e. 
$\Delta (d_i)=1\otimes d_i+d_i\otimes 1$.
\begin{lem}\label{l4.2} Relations $[X_{ik}+X_{jk},X_{ij}]=0$ for $i<j<k$, 
together with commutativity relations (\ref{3.2}), imply that the 
Jucys--Murphy elements  $d_j$ commute pairwise.
\end{lem}

{\it Proof.} For $j<l$,
$$[d_j,d_l]=\sum_{1\le i<j, \ 1\le k<l}[X_{ij},X_{kl}]=
\sum_{1\le i<j<l}\left\{ [X_{ij},X_{il}]+[X_{ij},X_{jl}]\right\} =0.
$$
\qed

Let us assume that additionally to the relations (\ref{3.1}) and 
(\ref{3.2}) the following relations are satisfied
\begin{equation}
X_{ij}^2=1,~~~1\le i<j\le n. \label{4.2}
\end{equation}
Then a map $(i<j)$
$$X_{ij}\to (i,j),
$$
where $(i,j)\in S_n$ is the transposition that interchanges $i$ and $j$ 
and fixes each $k\ne i,j$, defines a representation of the relations 
(\ref{3.1}), (\ref{3.2}) and (\ref{4.2}). In other words, there exists a 
homomorphism $p$ from the braid algebra $\B_n$ to the group ring of the 
symmetric group, $p:~\B_n\to\Z [S_n]$, such that $p(X_{ij})=(i,j)$.
Under this homomorphism $p$ the 
element $d_j$ maps to the Jucys--Murphy element in the group ring 
of the symmetric group $S_n$ (see, e.g. [J], 
[Mu], [Ra]). A proof follows from the following 
relations between the transpositions in the symmetric group $S_n$: if $1\le 
i<j<k\le n$, then
$$(ij)(ik)=(jk)(ij)=(ik)(jk),
$$
$$(ij)(jk)=(ik)(ij)=(jk)(ik).
$$

\begin{theorem}\label{t4.1} Let $\K_n$ be a commutative subalgebra of 
$\B_n$ generated by the Jucys--Murphy elements $d_j$, $j=2,\ldots ,n$. 
Then a map $d_j\mapsto x_{j-1}$ defines an isomorphism
$$\K_n\cong\Z [x_1,\ldots ,x_{n-1}].
$$
\end{theorem}

There exists a multiplicative analog, denoted by $D_k$, of the
Jucys--Murphy element $d_k$. 
Our construction of the elements $D_k$ follows to [Ra].

Let $g_i$, $1\le i\le n-1$, be the standard generators of the braid group 
$B_n$. Thus, the generators $g_i$ satisfy the following relations
\begin{eqnarray}
&& g_ig_j=g_jg_i, \ \ {\rm if} \ \ |i-j|\ge 2, \label{4.4'}\\
&&g_ig_{i+1}g_i=g_{i+1}g_ig_{i+1}, \ \ 1\le i\le n-2. \label{4.5'}
\end{eqnarray}
Now let us define the multiplicative Jucys--Murphy elements $D_k$: 
\begin{de} {\rm (cf. [Ra], (3.16))} The multiplicative Jucys--Murphy 
elements $D_k$ are defined by the following formulae
\begin{equation}
D_2=g_1^2; \ 
D_k=g_{k-1}g_{k-2}\cdots g_2g_1^2g_2\cdots g_{k-2}g_{k-1}, \ 3\le k\le n.
\label{4.6}
\end{equation}
\end{de}

\begin{lem}\label{l4.3} The elements $D_k$ commute pairwise.
\end{lem}

{\it Proof.} First of all, 
$$D_2D_3=g_1^2g_2g_1^2g_2=g_1g_2g_1g_2g_1g_2=g_2g_1^2g_2g_1^2=D_3D_2.
$$
Now, using induction, we have ($3\le k\le n-1$)
\begin{eqnarray*}
D_kD_{k+1}&=& g_{k-1}D_{k-1}g_{k-1}g_kg_{k-1}D_{k-1}g_{k-1}g_k=
g_{k-1}g_kD_{k-1}g_{k-1}D_{k-1}g_kg_{k-1}g_k\\
&=& g_{k-1}g_kD_{k-1}g_{k-1}D_{k-1}g_{k-1}g_kg_{k-1}=
g_{k-1}g_kD_{k-1}D_kg_kg_{k-1}\\
&=& g_{k-1}g_kD_kD_{k-1}g_kg_{k-1}=D_{k+1}D_k.
\end{eqnarray*}
Finally, if $l-k\ge 2$, then
\begin{eqnarray*}
D_kD_l&=& D_kg_{l-1}\cdots g_{k+1}D_{k+1}g_{k+1}\cdots g_{l-1}=
g_{l-1}\cdots g_{k+1}D_kD_{k+1}g_{k+1}\cdots g_{l-1}\\
&=& g_{l-1}\cdots g_{k+1}D_{k+1}D_kg_{k+1}\cdots g_{l-1}=D_lD_k.
\end{eqnarray*}
\qed

The elements $D_k$, $2\le k\le n$, 
generate a commutative subgroup $K_n$ in the braid group $B_n$. Now we 
are going to consider the following two problems: what is the 
infinitesimal deformation of the commutative subgroup $K_n$, and what is 
the quasi--classical limit
\begin{equation}
\Z [B_n]\to H_n(q)\to\Z [S_n] \label{4.3}
\end{equation}
of the multiplicative Jucys--Murphy elements $D_k$?  In the 
quasi--classical limit (\ref{4.3}) the 
group ring $\Z [B_n]$ of the braid group $B_n$ is degenerated at first to 
the Iwahori--Hecke algebra $H_n(q)$, and then to the group ring $\Z 
[S_n]$ of the symmetric group $S_n$.
\begin{pr}\label{p4.1} The quasi--classical limit (\ref{4.3}) of the 
multiplicative Jucys--Murphy element $D_k\in B_n$ is equal to the 
Jucys--Murphy element $d_k\in\Z [S_n]$.
\end{pr}

{\it Proof.}
Let us compute the quasi--classical limit (\ref{4.3}) of the element 
$D_k$. The first step is to consider $D_k$ as an element of the 
Iwahori--Hecke algebra $H_n(q)$. In other words, we have to add to 
(\ref{4.4'}) and (\ref{4.5'}) the new relations
\begin{equation}
g_i^2=(q^{\frac{1}{2}}-q^{-\frac{1}{2}})g_i+1, \ \ 1\le i\le n-1. 
\label{4.7}
\end{equation}
The next step is to consider the limit
$$\lim_{q\to 1}~\frac{D_k-1}{q^{\frac{1}{2}}-q^{-\frac{1}{2}}}.
$$
For this goal it is enough to compute $\ds\frac{dD_k}{dq}\vert_{q=1}$. 
This can be done using the formula
$$\frac{dg_i}{dq}=\frac{(q^{\frac{1}{2}}-q^{-\frac{1}{2}})g_i+
2}{2q(q^{\frac{1}{2}}+q^{-\frac{1}{2}})}.
$$
Thus,
\begin{eqnarray*}
\frac{dD_k}{dq}&=& \sum_{i=2}^{k-1}\left\{ g_{k-1}\cdots g_{i+1}
%\left(\frac{(q^{\frac{1}{2}}-q^{-\frac{1}{2}})g_i+
%2}{2(q^{\frac{1}{2}}+q^{-\frac{1}{2}})}\right) 
\left(\frac{dg_i}{dq}\right) g_{i-1}g_2g_1^2g_2\cdots g_{k-1}\right.\\ \\
&+&\left. g_{k-1}\cdots g_2g_1^2g_2\cdots g_{i-1}
%\left(\frac{(q^{\frac{1}{2}}-q^{-\frac{1}{2}})g_i
%+2}{2(q^{\frac{1}{2}}+q^{-\frac{1}{2}})}\right) 
\left(\frac{dg_i}{dq}\right) g_{i+1}\cdots g_{k-1}\right\} \\
&+& g_{k-1}\cdots g_2\left(\frac{d}{dq}g_1^2\right) g_2\cdots g_{k-1}\\ \\
&{\buildrel q=1\over=}& \sum_{i=2}^{k-1}\left\{\frac{1}{2}s_{k-1}
\cdots s_{i+1}s_is_{i+1}\cdots s_{k-1}+\frac{1}{2}s_{k-1}\cdots s_{i+1}
s_is_{i+1}\cdots s_{k-1}\right\}\\ \\
&+ & s_{k-1}\cdots s_2s_1s_2\cdots s_{k-1}\\ \\
&=& \sum_{i=1}^{k-1}(i,k)=d_k.
\end{eqnarray*}
\qed

Finally, let us explain a connection between the multiplicative 
Jucys--Murphy elements $D_k$ and the pure braid group $P_n$. More 
precisely, let us consider the infinitesimal deformations of the pure 
braid group $P_n$, and the multiplicative Jucys--Murphy elements $D_k$. 
To start, it is convenient to remind a few definitions.

By definition, the pure braid group $P_n$ is a kernel of the natural 
homomorphism
$$ B_n\to S_n, \ \ g_i\mapsto s_i=(i,i+1),
$$
where $s_i=(i,i+1)\in S_n$, $1\le i\le n-1$, denote the transposition 
that interchanges $i$ and $i+1$, and fixes all other elements of $[1,n]$.
We have an exact sequence:
$$1\to P_n\to B_n\to S_n\to 1.
$$

It is well--known (see, e.g. [B]), that $P_n$ is generated by the 
following elements
$$g_{ij}=(g_{j-1}g_{j-2}\cdots g_{i+1})g_i^2(g_{j-1}g_{j-2}\cdots 
g_{i+1})^{-1}, \ 1\le i<j\le n,
$$
subject the following relations: 
%where $1\le i<j\le n$, and $1\le k<m\le n$, and $j<m$:
\begin{eqnarray}
&& g_{ij}g_{kl}=g_{kl}g_{ij}, \  {\rm if~all} \ i,j,k,l \ {\rm 
are~distinct,} \label{4.8}\\
&& g_{ij}g_{ik}g_{jk}=g_{ik}g_{jk}g_{ij}=g_{jk}g_{ij}g_{ik}, \  {\rm if} \ 
1\le i<j<k\le n,\label{4.9}\\
&& g_{ik}g_{jk}g_{jl}g_{ij}=g_{jk}g_{jl}g_{ij}g_{ik}, \  {\rm if} \ 1\le 
i<j<k<l\le n. \label{4.10}
\end{eqnarray}
%\begin{equation}
%g_{ij}g_{km}g_{ij}^{-1}=\left\{\begin{array}{ll} g_{km},& k<i;\\ \\
%g_{jm}^{-1}g_{km}g_{jm}, & k=i;\\ \\
%g_{jm}^{-1}g_{im}^{-1}g_{jm}g_{im}g_{km}g_{im}^{-1}g_{jm}^{-1}
%g_{im}g_{jm}, & i<k<j;\\ \\
%g_{jm}^{-1}g_{im}^{-1}g_{km}g_{im}g_{jm}, & k=j;\\ \\
%g_{km} & k>j.
%\end{array}\right. 
%\end{equation}
One can show using only the relations (\ref{4.8})--(\ref{4.10}) that the 
elements
$$D_k=g_{1,k}g_{2,k}\cdots g_{k-1,k}\in P_n, \ 2\le k\le n,
$$ 
are mutually commute. 

Now let us consider the infinitesimal deformation 
$g_{ij}\mapsto 1+\epsilon X_{ij}$, of the pure braid group $P_n$. It is 
easy to see that the coefficients of $\epsilon^2$ on 
both sides of relations (\ref{4.8})--(\ref{4.10}) coincide with the 
defining relations (\ref{3.1})--(\ref{3.2}) 
for the braid algebra $\B_n$, and the element $D_k$ transforms to the 
Jucys--Murphy elements $d_k$ (more precisely, $D_k=1+\epsilon 
d_k+o(\epsilon^2)$).

\section{Quadratic Algebra ${\cal G}_n$}
\label{sec:quad}
\neweq

In the recent paper of the author and S.~Fomin [FK] the
quadratic algebra $\E_n^t$ was introduced. This algebra is closely related 
to the theory of 
quantum Schubert polynomials [FK], and their multiparameter deformation 
[K]. In this Section we introduce the quadratic algebra $\G_n$, which is a 
further generation of the quadratic algebra $\E_n^t$.

\begin{de}\label{d1.1} Define the algebra $\G_n$ (of type $A_{n-1}$) as the 
quadratic algebra (say, over $\Z$) with generators $[ij]$, $1\le i<j\le 
n$, which satisfy the following relations
\begin{equation} 
\label{1.1}
\begin{array}{rl}
&[ij][jk]=[jk][ik]+[ik][ij]\,,\\[.1in]
&[jk][ij]=[ik][jk]+[ij][ik],
                        \textrm{\ \ for } i<j<k\,;\quad\ \  \\[.1in]
\end{array}
\end{equation}
\begin{equation} 
\label{1.2}
\begin{array}{rl}
&~[ij][kl]=[kl][ij],\\[.1in]
&~\textrm{whenever } \{i,j\}\cap\{k,l\}=\phi\,,
   i<j, \textrm{ and }k<l \,. 
\end{array}
\end{equation}

\end{de}
The quadratic algebra $\E_n^t$, [FK], is the quotient of the algebra 
$\G_n$ by the two--side ideal generated by the relations
$$[ij]^2-t_{ij}=0,
$$
where the $t_{ij}$, for $1\le i<j\le n$, are a set of commuting 
parameters.

From (2.1) follows that the generators $[ij]$ satisfy the 
classical Yang--Baxter equation (CYBE)
\begin{equation}
[[ij],[ik]+[jk]]+[[ik],[jk]]=0, \ i<j<k, \label{2.2a}
\end{equation}
where the external brackets stand for the usual commutator: $[a,b]=ab-ba$.

\begin{theorem}\label{t1.1} The Hilbert series of the algebra $\G_n$ is 
given by
\begin{equation} 
\label{1.3}
H(\G_n;t)=\frac{1}{(1-t)(1-2t)\cdots (1-(n-1)t)}. 
\end{equation}
\end{theorem}
Theorem~\ref{t1.1} is a corollary of the following result which gives a 
description of an additive basis in the algebra $\G_n$:

\begin{theorem}\label{t1.2} A linear basis of the algebra $\G_n$ is given by 
the following set of noncommutative monomials
\begin{equation}
\label{1.4}
Z_2\cdot Z_3\cdot\ldots\cdot Z_n, 
\end{equation}
where
\begin{eqnarray*}
	Z_2 & = & \{\phi ,[12]\}  \\ 
	 & \cdots &  \\
	Z_k & = & \{\phi\}\cup\{ [i_1k]\cdots [i_lk]~\vert ~1\le i_j<k,~1\le 
	j\le l\} ,\ \ 2\le k\le n.
\end{eqnarray*}
\end{theorem}
A monomial in $[ij]$'s is called to be in {\it normal form}, if it 
belongs to the set (\ref{1.4}).

A proof of Theorem~\ref{t1.1} is similar to a proof of 
Corollary~\ref{c3.2} given by X.-S.~Lin, [L], Theorem~2.3. By 
this reason we omit a proof.

\begin{rem}\label{r1.2} {\rm It is clear 
that~~~~$H(Z_k;t)=\ds\frac{1}{1-(k-1)t}$, and }
\begin{equation}
H(\G_n;t)=\prod_{k=2}^nH(Z_k;t). \label{1.5}
\end{equation}
\end{rem}
\begin{cor}\label{c1.1} Let $\G_n'$ be a commutative version of the 
algebra $\G_n$, i.e. let us assume that the relations (\ref{1.2}) are 
valid for all $i,j,k,l$. Then
$$H(\G_n';t)=(1-t)^{-\left(\begin{array}{l}\small{n}\\\scriptsize{2}\end{array}\right)}.
$$
\end{cor}

\section{Dunkl elements}
\label{sec:dunkl}
\neweq

The Dunkl elements $\theta_j$, for $j=1,\ldots ,n$, in the quadratic 
algebra $\G_n$ are defined by
\begin{equation}
\label{2.1}
\theta_j=-\sum_{1\le i<j}[ij]+\sum_{j<k\le n}[jk]. 
\end{equation}
For the quadratic algebra $\E_n^t$ the Dunkl elements (\ref{2.1}) 
coincide with those introduced in [FK], Section~5.
The following result is well--known (cf., e.g. [Ki], Theorem~1.4, or 
[FK], Lemma~5.1)
\begin{lem}\label{l2.1} The classical Yang--Baxter equation (\ref{2.2a}), 
together 
with the commutativity relation (\ref{1.2}), imply that the Dunkl 
elements $\theta_j$ commute pairwise.
\end{lem}
There exists an obvious relation between Dunkl's elements $\theta_1,\ldots 
,\theta_n$ in the quadratic algebra $\G_n$, namely, $\theta_1+\cdots 
+\theta_n=0$. The result below shows that in the algebra $\G_n$ the Dunkl 
elements $\theta_1,\ldots ,\theta_{n-1}$ are algebraically independent 
(cf. Theorem~\ref{t5.1}).

\begin{theorem}\label{t2.1} Let $\H_n$ be a commutative subalgebra of $\G_n$ 
generated by the Dunkl elements $\theta_j$, $j=1,\ldots ,n$. Then a map 
$\theta_i\to x_i$, $1\le i\le n-1$, defines an isomorphism
$$\H_n\cong\Z [x_1,\ldots ,x_{n-1}],
$$
\end{theorem}
\begin{problem}\label{pr2.1} To describe the two--side ideals $\R\subset\G_n$ 
such that
\begin{equation}
H\left(\H_n/\R\cap\H_n;t\right)=[n]!:=\prod_{j=1}^n\frac{1-t^j}{1-t}.\label{2.2}
\end{equation}
\end{problem}

The examples of two--side ideals in the algebra $\G_n$ with the property 
(\ref{2.2}) will be given in Sections~7 and 8. These examples are related 
to the quantum cohomology ring of the flag variety [FK], and the 
multiparameter deformation of Schubert polynomials [K].

\begin{ex}\label{e2.1} {\rm Let us take $n=3$. In the algebra $\G_3$ we have 
the following relations
\begin{eqnarray*}
&&\theta_1+\theta_2+\theta_3=0,\\
&&\theta_1\theta_2+\theta_1\theta_3+\theta_2\theta_3+[12]^2+[13]^2+[23]^2=0,\\
&&\theta_1\theta_2\theta_3+[12]^2\theta_3+\theta_1[23]^2-[13]^2\theta_1-
\theta_3[13]^2=0,  \\
&&[\theta_2, [23]^2]=[\theta_1,[13]^2],
\end{eqnarray*}
where $[a,b]=ab-ba$ is the usual commutator.}
\end{ex}

\section{Dunkl and Jucys--Murphy elements}
\label{sec:djm}
\neweq

Let us start with the definition of the extended affine braid group 
$\wt B_n$.
\begin{de}\label{d6'.1} The extended affine braid group $\wt B_n$ is 
a group with generators
$$g_0, g_1,\ldots ,g_{n-1}, w
$$
which satisfy the following relations
\begin{eqnarray}
&& g_ig_j=g_jg_i, \ \   |i-j|\ge 2, \ 0\le i,j\le n-1, \label{6'.1}\\
&& g_ig_{i+1}g_i=g_{i+1}g_ig_{i+1}, \ \  0\le i\le n-1, \label{6'.2}\\
&& wg_i=g_{i-1}w,\ \ 0\le i\le n-1, \label{6'.3}
\end{eqnarray}
with indices understood as elements of $\Z /n\Z$.
\end{de}
It follows from (\ref{6'.3}) that $w^n$ is a central element.

There exists a canonical homomorphism $\pi$ from the extended affine 
braid group $\wt B_n$ to the classical braid group $B_n$. On the
generators $g_i$ and $w$ the homomorphism $\pi$ is given by the 
following rules
\begin{equation}
\pi (g_i)=g_i, \ 1\le i\le n-1, \label{6'.4}
\end{equation}
$$\pi (w)=g_{n-1}g_{n-2}\cdots g_2g_1.
$$
It follows from (\ref{6'.4}) that
$$\pi (g_0)=\pi (w)g_1\pi (w)^{-1}=g_{n-1}g_{n-2}\cdots 
g_2g_1g_2^{-1}\cdots g_{n-2}^{-1}g_{n-1}^{-1}.
$$

Now we are going to define the Dunkl 
ele\-ments $Y_1,\ldots ,Y_n$, and the dual Dunkl elements $Y_1^*,\ldots 
,Y_n^*$ in the extended affine braid group. 
\begin{de}\label{d6'.2} For each $i=1,\ldots ,n$, 
we define (cf. [KN]; [Ch1]) the Dunkl and dual Dunkl elements $Y_i$ and 
$Y_i^*$ respectively, by the following formulae
\begin{eqnarray}
Y_i& = & g_ig_{i+1}\ldots g_{n-1}wg_1^{-1}\ldots g_{i-1}^{-1}, 
\label{6'.5}\\
Y_i^*& = &g_i^{-1}g_{i+1}^{-1}\ldots g_{n-1}^{-1}wg_1\ldots g_{i-1}. 
\label{6'.6}
\end{eqnarray}
\end{de}

Note that $Y_1=g_1\ldots g_{n-1}w$, $Y_n=wg_1^{-1}\ldots g_{n-1}^{-1}$, 
$Y_1^*=wg_1\ldots g_{n-1}$ and\break $Y_n^*=g_1^{-1}\ldots g_{n-1}^{-1}w$.
\vskip 0.3cm

\hskip -0.6cm The Dunkl elements satisfy the following commutation relations with 
$g_1,\ldots ,g_{n-1}$:
\begin{eqnarray}
g_iY_{i+1}g_i &=& Y_i,~~~~~g_iY_j = Y_jg_i \ \ (j\ne i,i+1),
\label{6'.7}\\
g_iY_i^*g_i &=& Y_{i+1}^*, \ \ g_iY_j^* = Y_jg_i \ \ (j\ne i,i+1).\label{6'.8}
\end{eqnarray}
It is clear from our construction, that the image in the braid group 
$B_n$ of the dual Dunkl element $Y_k^*$, $2\le k\le n$, under the 
canonical homomorphism $\pi :\wt B_n\to B_n$ coincides with the 
multiplicative Jucys--Murphy element $D_k$, namely
\begin{eqnarray*}
\pi (Y_k^*) &= & g_k^{-1}g_{k+1}^{-1}\cdots g_{n-1}^{-1}\pi (w)g_1\cdots 
g_{k-1}\\
& = & g_{k-1}\cdots g_2g_1^2g_2\cdots g_{k-1}{\buildrel (\ref{4.6})\over 
=}D_k, \ \ 2\le k\le n,
\end{eqnarray*}
and $\pi (Y_1^*)=1$.
\begin{lem}\label{l6'.1} The Dunkl elements $Y_k$ (resp. the dual Dunkl 
elements $Y_k^*$) commute pairwise.
\end{lem}

Let us illustrate the main idea of the proof of Lemma~\ref{l6'.1} on some 
example. Let us take $n=5$ and prove that $Y_2^*Y_3^*=Y_3^*Y_2^*$. Indeed,
\begin{eqnarray*}
Y_2^*Y_3^* &=& g_2^{-1}g_3^{-1}g_4^{-1}wg_1g_3^{-1}g_4^{-1}wg_1g_2
=g_2^{-1}g_3^{-1}g_4^{-1}g_0g_2^{-1}g_3^{-1}g_4g_0w^2\\
& = & g_2^{-1}g_3^{-1}g_2^{-1}g_4g_3^{-1}g_0g_4g_0w^2
=g_3^{-1}g_2^{-1}g_3^{-1}g_4g_3^{-1}g_4g_0g_4w^2.
\end{eqnarray*}
Similarly,
$$Y_3^*Y_2^*=g_3^{-1}g_4^{-1}wg_1g_3^{-1}g_4^{-1}wg_1
=g_3^{-1}g_2^{-1}g_4^{-1}g_3^{-1}g_0g_4w^2.
$$
It is easy to see that in $\wt B_5$ we have
$$g_3^{-1}g_2^{-1}g_3^{-1}g_4g_3^{-1}g_4g_0g_4
=g_3^{-1}g_2^{-1}g_4^{-1}g_3^{-1}g_0g_4.
$$
\vskip 0.3cm

Now let us study the quasi--classical limit of the Dunkl
element $Y_k$. The first step is to consider $Y_k$ as an element of the 
extended affine Hecke algebra $H(\wt W)$, [KN], Section~2; [Ch1], [Ch2]. 
In other
words, we have to add to the relations (\ref{6'.1})--(\ref{6'.3}) the new
ones
\begin{equation}
g_i^2=(t-1)g_i+t, \ \ 0\le i\le n\label{6'.7'}
\end{equation}
where $t$ is a new parameter. The next step is to consider a 
representation of the extended affine Hecke algebra $H(\wt W)$ in the 
algebra $\D_{q,x}[W]$ of $q$--difference operators with permutations, 
see, for example, [Ch1], [KN]. Follow [KN], we define the elements 
$T_i$, $i=0,1,\ldots ,n-1$, in $\D_{q,x}[W]$ by
\begin{eqnarray*}
T_i&=&t+\frac{x_{i+1}-tx_i}{x_{i+1}-x_i}(s_i-1),\ \ i=1,\ldots ,n-1,\\
T_0&=&t+\frac{x_1-tqx_n}{x_1-qx_n}(s_0-1).
\end{eqnarray*}
Here $s_1,\ldots,s_{n-1}$ are the standard generators of the symmetric 
group $S_n$, and\break $s_0=ws_1w^{-1}$, where $w=s_{n-1}s_{n-2}\cdots s_1\tau_1$,
and $\tau_1=T_{q,x_1}$ is the $q$--shift operator:
$$(\tau_1f)(x_1,\ldots ,x_n)=f(qx_1,x_2,\ldots ,x_n).
$$
%it is clear that
%\begin{eqnarray*}
%(w^{-1}f)(x_1,\ldots ,x_n)&=&f(q_1x_n,x_1,x_2,\ldots ,x_{n-1}),\\
%(wf)(x_1,\ldots ,x_n)&=&f(x_2,x_3,\ldots , x_n,q^{-1}x_1),\\
%s_0f)(x_1,\ldots ,x_n)&=&f(qx_n,x_2,\ldots ,x_{n-1},q^{-1}x_1).
%\end{eqnarray*}

One can check that the elements $T_i$, $0\le i\le n-1$, and $w$, satisfy the
relations (\ref{6'.1})--(\ref{6'.3}) and (\ref{6'.7'}), and define
a representation of the extended affine Hecke algebra $H(\wt W)$ in the 
algebra $\D_{q,x}[W]$ of $q$--difference operators with permutations. 
In this representation, the Dunkl element $Y_k$ (resp. $Y_k^*$), 
$1\le k\le n$, up to a power of $t$ coincides (see [Ch1], [Ch2], [KN]) with 
the Dunkl--Cherednik operator ${\bf Y}_k$:
\begin{eqnarray}
{\bf Y}_k&=&t^{-n+2k-1}T_kT_{k+1}\ldots T_{n-1}wt_1^{-1}\ldots 
T_{k-1}^{-1}, \label{6'.8'}\\
{\bf Y}_k^*&=&t^{n-2k+1}T_k^{-1}T_{k+1}^{-1}\ldots T_{n-1}^{-1}w
T_1\ldots T_{k-1}. \label{6'.8"}
\end{eqnarray}
In order to understand better the quasi--classical behavior of the 
Dunkl--Cherednik operator ${\bf Y}_k$, it is convenient to rewrite 
slightly the formulae (\ref{6'.8'})--(\ref{6'.8"}). Namely, let us define 
($1\le i<j\le n$) the following operators acting on the ring of 
polynomials $\Z [t,t^{-1}][x_1,\ldots ,x_n]$
$$\overline T_{ij}=1+\frac{(1-t^{-1})x_j}{x_i-x_j}(1-s_{ij})=
t^{-1}T_{ij}s_{ij},
$$
where $s_{ij}$ is the exchange operator for variables $x_i$, $x_j$. Then
\begin{equation}
{\bf Y}_i=\overline T_{i,i+1}\overline T_{i,i+2}\ldots
\overline T_{i,n}\tau_i\overline T_{1,i}^{-1}\ldots
\overline T_{i-1,i}^{-1}, \label{6'.9'}
\end{equation}
where $\tau_i=s_is_{i+1}\ldots s_{n-1}ws_1\ldots s_{i-1}$ is the 
$q$--shift operator: $\tau_i=\tau_{x_i}=T_{q,x_i}$.
Finally, in the quasi--classical limit $q\to 1$ with 
rescaling $t=q^{\beta}$, we obtain from (\ref{6'.9'}) that
\begin{equation}
\D_j:=\lim_{q\to 1}~\frac{1-{\bf Y}_j}{1-q}=
x_j\frac{\partial}{\partial x_j}+\beta\sum_{i<j}x_j\partial_{ij}-
\beta\sum_{k>j}\partial_{jk}x_k. \label{6'.9"}
\end{equation}
Note that these Dunkl operators $\D_k$ commute with each other, i.e.
$[\D_i,\D_j]=0$.

More generally, let us define the affine extension of the algebra 
$\G_n$:
\begin{de} Define the algebra $\wt\G_n$ as the algebra (say over 
$\Z [q,q^{-1}]$)with generators
$$[ij], \ \ 0\le i<j\le n;\ \ x_i,\ 1\le i\le n, \ \ {\rm and} \ w,
$$
which satisfy the relations (\ref{1.1}) and (\ref{1.2}) for $i,j,k,l\in 
[1,n]$, and the 
following additional relations
\begin{eqnarray}
x_ix_j&=&x_jx_i, \ \ 1\le i,j\le n,\label{6'.9}\\
x_i[ab]&=&[ab]x_i, \ {\rm if} \ i\ne a,b, \label{6'.10}\\
x_j[ij]&=&[ij]x_i+1,\ \ x_i[ij]=[ij]x_j-1\ \ 1\le i<j\le n, \label{6'.11}\\
wx_i&=&q^{\delta_{1,i}}x_{i-1}w, \label{6'.11'}\\
w[ij]&=&q^{-\delta_{i,0}}[i-1,j-1]w, \label{6'.11''}
\end{eqnarray}
with indices understood as elements of $\Z /n\Z$.
\end{de}
It is clear from (\ref{6'.11}), that
\begin{equation}
[x_i+x_j,[ij]]=0. \label{6'.12}
\end{equation}
We define the Dunkl elements $\wt\theta_j$ in the algebra $\wt\G_n$ by the 
following rule
\begin{equation}
\wt\theta_j=x_j+\theta_j, \ \ 1\le j\le n. \label{6'.13}
\end{equation}
The following result is well-known, see, e.g. [Ki]:
\begin{lem}\label{l6'.2} The classical Yang--Baxter equation (\ref{2.2a}), 
together with commutativity relation (\ref{1.2}) and relation 
(\ref{6'.12}), imply that the Dunkl elements $\wt\theta_j=x_j+\theta_j$ 
commute pairwise.
\end{lem}

{\it Proof.} If $i<j$, then
$$[\wt\theta_i,\wt\theta_j]=[x_i,x_j]+[x_i,\theta_j]+[\theta_i,x_j]+
[\theta_i,\theta_j]=[x_i+x_j,[ij]]=0.
$$
\qed

Let us introduce the algebra $\wt\G_n^0$ which is the quotient of the 
algebra $\wt\G_n$ by the two--side ideal generated by the relation
\begin{equation}
[ij]^2=0, \ \ 1\le i<j\le n. \label{6'.14}
\end{equation}

\begin{de}\label{d6'.3} For $1\le i<j\le n$ let us define an element 
$T_{ij}\in\wt\G_n^0$ by the following formula (cf. [KN])
\begin{equation}
T_{ij}=t-(x_j-tx_i)[ij]. \label{6'.15}
\end{equation}
\end{de}

It follows from (\ref{6'.11}) and (\ref{6'.14}), that
$$T_{ij}^2=(t-1)T_{ij}+t.
$$
Our next goal is to check that the elements 
$\wt T_i:=T_{i,i+1}$, $1\le i\le n-1$, satisfy the Coxeter relations.
\begin{pr}\label{p6'.1} We have the following relation in the algebra 
$\wt\G_n^0$:
\begin{equation}
T_{ab}T_{bc}T_{ab}=T_{bc}T_{ab}T_{bc}, \ \ a<b<c. \label{6'.16}
\end{equation}
\end{pr}

{\it Proof.} Direct computation based on (\ref{6'.11}) and the following 
relation in the algebra $\wt\G_n^0$:
$$[ab][bc][ab]=[bc][ab][bc], \ \ 1\le a<b<c\le n.
$$
\qed

As a corollary, we see that the elements $w$, $\wt T_0=w\wt T_1w^{-1}$, 
$\wt T_1,\ldots ,\wt T_{n-1}$ generate 
the representation of the extended affine Hecke algebra $H(\wt W)$. In 
this representation the quasi--classical limit $q\to 1$ with rescaling 
$t=q^{\beta}$ of the Dunkl element $t^{-n+2j-1}Y_j$ is equal to
\begin{equation}
x_j\frac{\partial}{\partial x_j}+\beta\sum_{i<j}x_j[ij]-
\beta\sum_{k>j}[jk]x_k. \label{6'.17}
\end{equation}
The element (\ref{6'.17}) does not belong to the algebra $\wt\G_n$, but 
its further extension $\wt{\wt\G}_n$ which is an analog of the double 
affine Hecke algebra, introduced by I.~Cherednik. It is an interesting 
problem to understand some connections between the Schubert calculus and 
the Dunkl--Cherednik operators.

\begin{con} {\rm (Nonnegativity conjecture, cf. [FK], Conjecture~8.1)} 
For any $w\in S_n$, the Schubert polynomial ${\s}_w$ evaluated at the Dunkl 
elements $\wt\theta_1,\ldots ,\wt\theta_n$ belongs to the positive cone 
$\wt\G_n^+$, where $\wt\G_n^+$ is the cone of all nonnegative integer 
linear combinations of all (noncommutative in general) monomials in the 
generators $x_i$, $1\le i\le n$, and $[ij]$, $1\le i<j\le n$, of the 
algebra $\wt\G_n^0$.
\end{con}

\newpage
\section{Quadratic algebra $\E_n^t$}
\label{sec:ent}
\neweq

In this and the next sections we study some interesting quotients of the the 
quadratic algebra $\G_n$.

\begin{de}\label{d5.1} {\rm (cf. [FK], Section~15)} Define the algebra 
$\E_n^t$ as 
the algebra over $\Z$ with generators $[ij]$, $1\le i<j\le n$, which 
satisfy the relations (\ref{1.1}), (\ref{1.2}) and the additional relations
\begin{equation}
([ij],[kl]^2)=0, ~~for~ all~~i,j,k,l.
\label{5.1}
\end{equation}
\end{de}
The relations (\ref{5.1}) mean that the squares $[kl]^2$ belong to the 
center of $\E_n^t$. Let us put $t_{ij}=[ij]^2$, and consider the ring of 
polynomials $\Z [t]:=\Z [~t_{ij},~1\le i<j\le n]$.

\begin{con}\label{c5.1} {\rm ([FK])} The algebra 
$\E_n^t$ is a finite dimensional module over the ring of polynomials $\Z 
[t]$.
\end{con}

\begin{theorem}\label{t5.1} {\rm ([FK], Conjecture~15.1; [P])} Let 
$\H_n^t$ be a commutative subalgebra of $\E_n^t$ generated by the Dunkl 
elements $\theta_1,\ldots ,\theta_n$. Then
$$\H_n^t\cong\Z [t][\theta_1,\ldots ,\theta_n]/J,
$$
where the ideal $J$ is generated by the generalized elementary functions
\begin{equation}
e_m(X_n|t):=\sum_l\sum_{\begin{array}{c} \scriptstyle{1\le i_1<\cdots 
<i_l\le n}\\ \scriptstyle{j_1>i_1,\ldots ,j_l>i_l}\end{array}}
e_{m-2l}(X_{\overline{I\cap J}})\prod_{k=1}^lt_{i_kj_k},
\label{5.2}
\end{equation}
where $i_1,\ldots ,i_l,j_1,\ldots ,j_l$ should be distinct elements of the
set $\{ 1,\ldots ,n\}$, and $X_{\overline{I\cap J}}$ denotes set of 
variables $x_a$, for which the subscript $a$ is neither one of the $I_k$
nor one of the $J_k$.
\end{theorem}
One can check that the relations (\ref{5.1}) are 
equivalent to the following ones
$$(\theta_i,[kl]^2)=0, \ \ {\rm for~ all} \ \ i,k,l.
$$
The relations (\ref{5.1}) in the algebra $\E_n^t$ imply some nontrivial 
relations between monomials in the space 
$Z:=Z_2\cdot Z_3\cdot\ldots\cdot Z_n$. We will 
describe part of these relations in Proposition~\ref{p5.1} below, but before 
that, let us introduce some additional notations. For given $n\ge 3$ and 
$k$, $1\le k<n$, let us denote by $\Lambda_{n,k}$ a set of all sequences 
of integer numbers $A=(a_1,\ldots ,a_k)$, such that $a_1+\cdots +a_k=n-2$, 
and $a_i\in\Z_{>0}$. For each sequence $A$ from the set $\Lambda_{n,k}$ 
let us define a monomial $[A]$ in the space $Z$ as the following ordered 
product
$$[A]=\prod_{j=1}^k\prod_{l=1}^{a_j}[n-k+l-2-\sum_{s=1}^j(a_s-1), 
~n-k+j-1].
$$
For example, for $n=7$ and $A=(2,1,2)$ we have
$$[A]=[24][34][25][16][26].
$$
\begin{pr}\label{p5.1}
\begin{equation}
F_n:=\sum_{i=1}^{n-1}[i,n][i+1,n]\cdots [n-1,n][1,n][2,n]\cdots [i,n]
\label{5.3}
\end{equation}
$$=\sum_{k=2}^{n-1}(-1)^{n-k-1}(t_{kn}-t_{k-1,n})
\sum_{A\in\Lambda_{n,n-k}}[A].~~~~~
$$
\end{pr}

{\it Proof.} By induction, using the following 
formula 
\begin{eqnarray*}
&& \Phi_{n+1}\cdot [12]-[12]\cdot\Phi_{n+1}=F_{n+1},\ \ {\rm where}\\
&&\Phi_n:=[2,n]\cdots [n-1,n][2,n]+\cdots +[n-1,n][2,n]\cdots 
[n-2,n][n-1,n].
\end{eqnarray*}

\qed

Let us make the following simple but useful observation. The action of 
the symmetric group, namely, $w([i_1j_1]\cdots [i_kj_k])=[w(i_1)w(j_1)]
\cdots [w(i_k)w(j_k)]$, $w\in S_n$, transforms every specific identity in 
$\G_n$ (resp. $\E_n^t$, $\E_n^0$) into an orbit of identities. As an 
example, the identity (in $\E_n^t$)
$$[12][13][14][12]+[13][14][12][13]+[14][12][13][14]=
(t_{34}-t_{24})[13][23]-(t_{24}-t_{14})[12][13]
$$
(a special case of Proposition~\ref{p5.1}) immediately implies a more general
identity
\begin{eqnarray*}
&&[ab][ac][ad][ab]+[ac][ad][ab][ac]+[ad][ab][ac][ad]\\
&&=([cd]^2-[bd]^2)[ac][bc]-([bd]^2-[ad]^2)[ab][ac].
\end{eqnarray*}

Proposition below describes some interesting relations in the quadratic 
algebra $\E_n^t$. They can be interpreted as the Pieri rules for the 
classical and quantum Schubert polynomials, [FK], [P], and for the 
multiparameter deformation of Schubert polynomials [K], [P].

\begin{pr}\label{p5.2} {\rm ([P])} Let $A$ be a subset in $[1,\ldots ,n]$, 
and let us put $\theta_A=\ds\prod_{i\in A}\theta_i$. Then in the algebra 
$\E_n^t$ we have the following relation
$$e_m(\theta_A|t)=\sum [i_1j_1]\cdots [i_mj_m],
$$
where the summation is taken over all sequences $I=\{ i_1,\ldots ,i_m\}$,
$J=\{ j_1,\ldots ,j_m\}$ such that

i)~~~ $I\subset A$, \ $J\subset [1,n]\setminus A$,

ii)~~ $i_1,\ldots ,i_m$ are all distinct,

iii)~  $1\le j_1\le j_2\le\cdots\le j_m\le n$.
\end{pr}

\begin{problem} Find a nontrivial faithful representation of the 
quadratic algebra $\E_n^t$.
\end{problem}

\section{Quadratic algebra $\E_n^0$ }
\label{sec: En0}
\neweq

In this Section we study the properties of the algebra $\E_n^t$ when all 
parameters $t$ are equal to zero. Hence, in this Section we assume that
\begin{equation}
[ij]^2=0, \ \ 1\le i<j\le n.
\label{6.1}
\end{equation}
This algebra was introduced and studied at first in [FK]. It follows 
from Theorem~\ref{t1.2} that any element of $\E_n^0$ can be presented as 
$\Z$-linear combination of monomials $z$ of the following form
\begin{equation}
z=z_1z_2\cdots z_n,
\label{6.2}
\end{equation}
where $z_i\in Z_i$, and ($2\le k\le n$)
\begin{equation}
Z_k:=Z_{n,k}=\{\phi\}\cup\{ [i_1k]\cdots [i_lk]~|~1\le i_l<k, \ 1\le 
j\le l\}. \label{6.3}
\end{equation}

In the algebra $\G_n$ the monomials in the set $Z_k$ are linearly 
independent (Theorem~\ref{t1.2}). However, in the quotient algebra $\E_n^0$
%(or $\E_n^t$) 
some new relations between the monomials in $Z_k$ are appeared. The 
Proposition~\ref{p6.1} below is a special case of Proposition~\ref{p5.1} and
gives the relations in subalgebra $Z_m\subset\E_n^0$, but not only
in the space $Z=Z_2\cdot Z_3\cdot\ldots\cdot Z_n$.

\begin{pr}\label{p6.1} {\rm (cf. [FK])} Let $a_1,\ldots ,a_m$ be a sequence 
of pairwise distinct integer numbers, $1\le a_i<n$, $m\le n$. Then
\begin{equation}
a_1a_2\cdots a_ma_1+a_2a_3\cdots a_1a_2+\cdots +a_ma_1a_2\cdots 
a_{m-1}a_m=0.
\label{6.4}
\end{equation}
Here we used an abbreviation $a_j=[a_jn]$.
\end{pr}

{\it Proof.} It is enough to prove that
\begin{equation}
F_n:=\sum_{i=1}^{n-1}[i,n][i+1,n]\cdots [n-1,n][1,n][2,n]\cdots [i,n]=0.
%12\cdots n-11+23\cdots n-112+\cdots +n-112\cdots n-2n-1=0.
\label{6.5a}
\end{equation}
Let us prove (\ref{6.5a}) by induction. Thus we may assume that $F_k=0$ 
for $k\le n$. Let us consider the following element in $Z_{n+1}$:
$$\Phi_n=[2,n+1]\cdots [n,n+1][2,n+1]+\cdots +[n,n+1][2,n+1]\cdots 
[n-1,n+1][n,n+1].
$$

By induction assumption, we have $\Phi_n=0$. On the other hand, it is 
easy to check that
\begin{equation}
\Phi_n\cdot [12]-[12]\cdot\Phi_n=F_{n+1}. \label{6.5b}
\end{equation}
The last equality shows that $F_{n+1}=0$.
\qed

\begin{ex} {\rm If $m=3$, the relation (\ref{6.4}) has the following form 
($a<b<n$)
\begin{equation}
aba+bab=0, \ \ {\rm in~the~algebra}\ \ \E_n^0
\label{6.5}
\end{equation}
Similarly, if $a,b,c$ are distinct, $<n$, then}
\begin{equation}
abca+bcab+cabc=0, \ \ {\rm in~the~algebra}\ \ \E_n^0
\label{6.6}
\end{equation}
\end{ex}
\vskip 0.3cm

We don't know how to describe all
relations between monomials from the space $Z_k$ for fixed $k$, but 
we can show that if 
$k\le 5$, then all relations between monomials in $Z_k$ follows from 
(\ref{6.4}). This observation allows us to compute the dimensions of the 
homogeneous components $\E_{n,k}^0$ for $1\le k\le 5$, and also 
$\dim\E_{n,6}^0$. We summarize the 
results in the following
\begin{pr}\label{p6.2}~

\vskip 0.3cm

$\bullet$ ${\rm dim}\,\E_{n,1}^0=\left(\begin{array}{c}n\\ 
2\end{array}\right)$,
\vskip 0.3cm

$\bullet$ ${\rm dim}\,\E_{n,2}^0=3\left(\begin{array}{c}n\\ 
4\end{array}\right) +4\left(\begin{array}{c}n\\ 3\end{array}\right)$,
\vskip 0.3cm

$\bullet$ ${\rm dim}\,\E_{n,3}^0=15\left(\begin{array}{c}n\\ 
6\end{array}\right) +40\left(\begin{array}{c}n\\ 5\end{array}\right) +
30\left(\begin{array}{c}n\\ 4\end{array}\right) +3\left(\begin{array}{c}n\\ 
3\end{array}\right)$,
\vskip 0.3cm

$\bullet$ ${\rm dim}\,\E_{n,4}^0=105\left(\begin{array}{c}n\\ 
8\end{array}\right) +420\left(\begin{array}{c}n\\ 7\end{array}\right) +
610\left(\begin{array}{c}n\\ 6\end{array}\right) +366\left(\begin{array}{c}n\\ 
5\end{array}\right)$
\vskip 0.3cm

\hskip 2cm $+67\left(\begin{array}{c}n\\ 
4\end{array}\right) +\left(\begin{array}{c}n\\ 
3\end{array}\right)$,
\vskip 0.3cm

$\bullet$ ${\rm dim}\,\E_{n,5}^0=945\left(\begin{array}{c}n\\ 
10\end{array}\right) +5040\left(\begin{array}{c}n\\ 9\end{array}\right) +
10780\left(\begin{array}{c}n\\ 8\end{array}\right) +
11571\left(\begin{array}{c}n\\ 7\end{array}\right)$
\vskip 0.3cm

\hskip 2cm $+6285\left(\begin{array}{c}n\\ 6\end{array}\right) +
1480\left(\begin{array}{c}n\\ 5\end{array}\right) +
96\left(\begin{array}{c}n\\ 4\end{array}\right)$,
\vskip 0.3cm

$\bullet$ ${\rm dim}\,\E_{n,6}^0=10395\left(\begin{array}{c}n\\ 
12\end{array}\right) +69300\left(\begin{array}{c}n\\ 11\end{array}\right) +
195300\left(\begin{array}{c}n\\ 10\end{array}\right) +
299908\left(\begin{array}{c}n\\ 9\end{array}\right)$
\vskip 0.3cm

\hskip 2cm $+268674\left(\begin{array}{c}n\\ 8\end{array}\right) +
138545\left(\begin{array}{c}n\\ 7\end{array}\right) +
37456\left(\begin{array}{c}n\\ 6\end{array}\right) +
4231\left(\begin{array}{c}n\\ 5\end{array}\right)$
\vskip 0.3cm

\hskip 2cm $+106\left(\begin{array}{c}n\\ 4\end{array}\right)$.
\end{pr}
\vskip 0.3cm
\begin{con}\label{c6.1} 
\begin{equation}
{\rm dim}\,\E_{n,k}^0=(2k-1)!!\left(\begin{array}{c}n\\ 2k\end{array}\right) +
(2k-1)!!\frac{4(k-1)}{3}\left(\begin{array}{c}n\\ 2k-1\end{array}\right)
\label{6.7}
\end{equation}
$$~~~~~~~~~~~~~~+(2k-3)!!\frac{(k-1)(k-2)(16k-3)}{9}
\left(\begin{array}{c}n\\ 2k-2\end{array}\right) +\cdots .
$$
\end{con}

{\it Proof of Proposition~\ref{p6.2}}. First of all, let us compute ${\rm 
dim}\,Z_{n,k}$ for small values of $k$. It is clear that ${\rm 
dim}\,Z_{n,1}=n-1$ and ${\rm dim}\,Z_{n,2}=(n-1)(n-2)$, because there are 
no nontrivial relations of degree $<3$ in $Z_n$. In degree 3 we have to 
take into account the relations (\ref{6.5}). Thus,

\vskip 0.3cm
${\rm dim}\, Z_{n,3}=(n-2){\rm dim}\, Z_{n,2}-\ds\left(\begin{array}{c}n-1\\ 
2\end{array}\right) =\frac{(n-1)(n-2)(2n-5)}{2}.$
\vskip 0.3cm

Similarly, in degree 4 we have to consider both the relations (\ref{6.5}) 
and (\ref{6.6}). Thus,

${\rm dim}\, Z_{n,4}=(n-2){\rm dim}\, Z_{n,3}-(n-2)\left(\begin{array}{c}n-1\\ 
2\end{array}\right) -2\left(\begin{array}{c}n-1\\ 
3\end{array}\right) $

\vskip 0.3cm
\hskip 1.6cm$=\ds\frac{(n-1)(n-2)(n-3)(3n-7)}{3}=24\left(\begin{array}{c}
n-1\\ 4\end{array}\right) +10\left(\begin{array}{c} n-1\\ 
3\end{array}\right) .$
\vskip 0.3cm

Using the same method, we can compute
\vskip 0.3cm

${\rm dim}\, Z_{n,5}=120\left(\begin{array}{c}n-1\\ 5\end{array}\right) +
86\left(\begin{array}{c}n-1\\ 4\end{array}\right) +
9\left(\begin{array}{c}n-1\\ 3\end{array}\right) . $

\vskip 0.3cm
When $k\ge 6$, some additional to (\ref{6.4}) relations between monomials 
in $Z_k$ will appear. We describe such new relations in degree 6.
\begin{lem}\label{l6.3} Let $a,b,c,d$ be a set of distinct integer 
numbers. Then 
\begin{eqnarray}
&& abcdca+bcdcab+cabadc+dcabad+cdcaba \label{8.9"}\\
&& =abacdc+acdcba+bacdcb+cdabac+dabacd.\nonumber
\end{eqnarray}
where $a:=[an]$, \ $n\ge 5$.
\end{lem}

{\it Proof of Lemma~\ref{l6.3}.} Let us consider the difference
\begin{equation}
[ac]F_4(a,d,c,b)-F_4(a,b,c,d)[ac], \label{8.9'}
\end{equation}
where
\begin{eqnarray*}
F_4(a,b,c,d)&=&[an][bn][cn][dn][an]+[bn][cn][dn][an][bn]\\
&+&[cn][dn][an][bn][cn]+[dn][an][bn][cn][dn].
\end{eqnarray*}
First of all, the difference (\ref{8.9'}) is equal to zero. On the other 
hand, using the relations in the algebra $\E_n^0$, one can transform the 
expression (\ref{8.9'}) to the difference between the left and 
the right hand sides of (\ref{8.9"}).
\qed

Using Lemma~\ref{l6.3}, we can compute

\vskip 0.4cm
${\rm dim}\, Z_{n,6}=720\left(\begin{array}{c}n-1\\ 6\end{array}\right) +
756\left(\begin{array}{c}n-1\\ 5\end{array}\right) +
187\left(\begin{array}{c}n-1\\ 4\end{array}\right) +
6\left(\begin{array}{c}n-1\\ 3\end{array}\right) .$
\vskip 0.2cm

Finally, in order to compute the dimensions of homogeneous components 
$\E_{n,k}^0$, $1\le k\le 6$, we use an observation that monomials from 
the set $Z_2\cdot Z_3\cdot Z_4\cdot Z_5\cdot Z_6$ form a basis in 
$\displaystyle\bigoplus_{k=0}^6\E_{n,k}^0$. 
The proof of Proposition~\ref{p6.2} is 
finished.
\qed
\begin{con}\label{c6.2} {\rm (Factorization property)}
\begin{equation}
H(\E_n^0;t)=H(Z_{n,2};t)H(Z_{n,3};t)\cdots H(Z_{n,n};t). \label{6.8}
\end{equation}
\end{con}
\begin{ex}\label{e6.1} {\rm If $n=2$, then $Z_2=\{\phi\}\cup\{ 12\}$, and 
$H(Z_2;t)=1+t$.

\vskip 0.3cm
If $n=3$, then
$Z_3=\{\phi\}\cup\{ 13, \ 23, \ 13\cdot 23,\  23\cdot 13, \ 
13\cdot 23\cdot 13\},$ and
$$H(Z_3;t)=1+2t+2t^2+t^3.
$$

Now, let us consider $n=4$. Then ($a:=[an]$)
\begin{eqnarray*}
Z_{4,1}& = & \{ 1, 2, 3\},\\
Z_{4,2}& = & \{ 12, 13, 21, 31, 23, 32\},\\
Z_{4,3}& = & \{ 121, 123, 131, 132, 213, 312, 231, 232, 321\},\\
Z_{4,4}& = & \{ 1213, 1231, 1232, 1312, 1321, 2131, 2132, 3121, 2312, 2321\},\\
Z_{4,5}& = & \{ 12131, 12132, 12312, 12321, 13121, 21312, 21321, 31213, 
23121\},\\
Z_{4,6}& = & \{ 121312, 121321, 123121, 131213, 213121, 231213\},\\
Z_{4,7}& = & \{ 1213121, 1231213, 2131213\},\\
Z_{4,8}& = & \{ 12131213\}.
\end{eqnarray*}
Hence, $H(Z_4;t)=(1+t)^2(1+t+t^2)(1+t^2)^2$, and 
$$H(\E_4^0;t)=(1+t)^4(1+t^2)^2(1+t+t^2)^2.
$$

Formula for $H(\E_4^0;t)$ was appeared in [FK] and confirmed by J.-E.~Ross
[Ro], who also computed the Hilbert polynomial $H(\E_5^0;t)$.}
\end{ex}

\section{Quadratic algebra $\B_n^t$ }
\label{sec: Bnt}
\neweq

In this Section we define a certain quotient $\B_n^t$ of the braid 
algebra $\B_n$, and formulate a conjecture that the algebras $\B_n^t$ and 
$\E_n^t$ have the same Hilbert polynomials.

\begin{de}\label{d10.1} define the algebra $\B_n^t$ as the algebra over 
$\Q$ with generators $X_{ij}$, $1\le i<j\le n$, which satisfy the 
relations (\ref{3.1}), (\ref{3.2}) and the additional relations
\begin{equation}
[ij]^2=t_{ij}, \label{10.1}
\end{equation}
where the $t_{ij}$ for $1\le i<j\le n$, are a set of commuting parameters.
\end{de}
In the case when all parameters $t_{ij}$ are equal to zero, we denote the 
algebra $\B_n^t$ by $\B_n^0$.

\begin{con}\label{c10.1} The algebra $\B_n^t$ is a finite dimensional 
module over the ring of polynomials $\Q [t]:=\Q [t_{ij}, 1\le i<j\le n]$.
\end{con}

\begin{problem} To find all relations among the Jucys--Murphy 
elements $d_j$,\break $2\le j\le n$, and find the Hilbert polynomial of 
commutative subalgebra $\K_n^0$ of the algebra $\B_n^0$ generated by the 
Jucys--Murphy elements $d_j$, $2\le j\le n$.
\end{problem}

\begin{ex}\label{e10.1} {\rm Let us consider the algebra $\B_3^0$. 
This is an algebra over $\Q$ with generators $X_{12},X_{13},X_{23}$ 
subject the following quadratic relations
\begin{eqnarray}
&& X_{12}^2=X_{13}^2=X_{23}^2=0, \label{10.2}\\
&& X_{13}X_{12}=X_{12}X_{13}-X_{13}X_{23}+X_{23}X_{13}, \label{10.3}\\
&& X_{23}X_{12}=X_{12}X_{23}-X_{23}X_{13}+X_{13}X_{23}.~~~~~~~~~~~~~
~~~~~~~~~~~~\label{10.4}
\end{eqnarray}
Using these relations, we can find the new relations in the algebra 
$\B_3^0$, namely,
\begin{eqnarray}
&& X_{23}X_{13}X_{23}-X_{13}X_{23}X_{13}+X_{12}X_{13}X_{23}-
X_{12}X_{23}X_{13}=0, \label{10.5}\\
&& X_{23}X_{13}X_{23}X_{13}+X_{12}X_{13}X_{23}X_{13}=0, \label{10.6}\\
&& X_{13}X_{23}X_{13}X_{23}+X_{12}X_{23}X_{13}X_{23}=0, \label{10.7}\\
&& X_{12}X_{23}X_{13}X_{23}-X_{12}X_{13}X_{23}X_{13}=0. \label{10.8}
\end{eqnarray}
\begin{rem}\label{r10.1} {\rm To be more precise, from the relations 
(\ref{10.2})--(\ref{10.4}) follow only that
$$2\cdot{\rm LHS}(9.5)=0
$$
We expect that the $\Z$--form of the quadratic algebra $\B_n^0$ has only 
2--torsion.}
\end{rem}

It is easy to check from the relations (\ref{10.2})--(\ref{10.8}) that 
the Hilbert polynomial of the algebra $\B_3^0$ is equal to
$$H(\B_3^0;t)=1+3t+4t^2+3t^3+t^4=(1+t)^2(1+t+t^2)=H(\E_3^0;t).
$$
Now, let us consider the commutative subalgebra $\K_3^0:=\Q [d_2,d_3]$, 
where $d_2=X_{12}$ and $d_3=X_{13}+X_{23}$. The Jucys--Murphy elements 
$d_2$ and $d_3$ satisfy the following relations
$$d_2^2=0,\ \ d_2d_3=d_3d_2, \ \ (d_2+d_3)d_3^3=0,
$$
and the elements $1,d_2,d_3,d_2d_3,d_3^2,d_2d_3^2,d_3^3,d_2d_3^3$ form a 
linear basis in the algebra $\K_3^0$. Thus we have
$$\dim\K_3^0=8, \ \ H(\K_3^0;t)=(1+t)^2(1+t^2).
$$
For general $n$ we expect that
\begin{equation}
H(\K_n^0;t)=\prod_{k=1}^{n-1}\frac{1-t^{2k}}{1-t}=H(W(B_{n-1});t),
\end{equation}
where $W(B_{n-1})$ is the Weyl group of type $B_{n-1}$.}
\end{ex}

\begin{con}\label{c10.2} The quadratic algebras $\B_n^0$ and $\E_n^0$ have 
the same Hilbert polynomials
\begin{equation}
H(\B_n^0;t)=H(\E_n^0;t). \label{10.9}
\end{equation}
\end{con}
We can check (\ref{10.9}) for $n=3$ (Example~\ref{e10.1}) and $n=4$. In 
order to prove the equality (\ref{10.9}) for $n=4$, we use the same method 
as in Example~\ref{e6.1}. But, instead of identity (\ref{6.6})  in the 
algebra $\E_n^0$, we have to use the following 14--terms relations in the 
algebra $\B_n^0$, $n\ge 4$:
\begin{eqnarray}
&& X_{23}(X_{14}X_{24}X_{34}-X_{24}X_{14}X_{34}-X_{34}X_{14}X_{24}+
X_{34}X_{24}X_{14}) \label{10.10}\\
&&+(X_{13}X_{23}-X_{23}X_{13})(X_{24}X_{34}-X_{34}X_{24})
-X_{14}X_{24}X_{34}X_{24} \nonumber\\
&&+X_{24}X_{34}X_{24}X_{14}+X_{24}X_{14}X_{34}X_{24}
-X_{24}X_{34}X_{14}X_{24}\nonumber\\
&&+X_{34}X_{14}X_{24}X_{34}-X_{34}X_{24}X_{14}X_{34}=0, \nonumber \\ 
\nonumber\\
&&X_{12}(X_{14}X_{34}X_{24}-X_{24}X_{14}X_{34}-X_{34}X_{14}X_{24}+
X_{24}X_{34}X_{14}) \label{10.11}\\
&&+(X_{13}X_{23}-X_{23}X_{13})(X_{14}X_{24}-X_{24}X_{14})
+X_{14}X_{24}X_{34}X_{14} \nonumber\\
&&-X_{14}X_{34}X_{24}X_{14}-X_{14}X_{24}X_{14}X_{34}
+X_{24}X_{14}X_{34}X_{24}\nonumber\\
&&-X_{24}X_{34}X_{14}X_{24}+X_{34}X_{14}X_{24}X_{14}=0. \nonumber
\end{eqnarray}

It is also natural to ask whether or not the quadratic algebra $\B_n^0$ 
has the Koszul property (see, e.g. [Ma] and references therein).

%\newpage
\section{Quadratic algebra $\A_n^t$ }
\label{sec: quadr alg}
\neweq

In this section we consider the quadratic algebra $\A_n^t$ (denoted by
$\E C_n$ in [FK], Section~4.3) which is the 
commutative quotient of the algebra $\E_n^t$.

\begin{de}\label{d7.1} Define the algebra $\A_n^t$ (of type $A_{n-1}$) as 
the quadratic algebra over $\Z$ with generators $[ij]$, $1\le i<j\le n$, 
which satisfy the following relations
\begin{eqnarray}
&&[ij][kl]=[kl][ij], \ \ {\rm for~ all} \ i,j,k,l; \label{7.1}\\
&&[ij][jk]=[ik][ij]+[ik][jk], \ \ i<j<k; \label{7.2}\\
&&[ij]^2=t_{ij}. \label{7.3}
\end{eqnarray}
\end{de}
Let us denote by $\A_n^0$ 
the specialization $t_{ij}=0$, $1\le i<j\le n$, of the algebra $\A_n^t$.

\begin{theorem}\label{t7.1} {\rm (cf. [A])} The Hilbert polynomial of the 
algebra $\A_n^0$ is given by
\begin{equation}
H(\A_n^0;t)=(1+t)(1+2t)\cdots (1+(n-1)t). \label{7.4}
\end{equation}
\end{theorem}
 
\begin{cor}\label{c7.1} {\rm (cf. [A])} An additive basis of the 
algebra $\A_n^t$ is given by the following set of (commutative) monomials
$$[j_1k_1]\ldots [j_lk_l],
$$
where $j_s<k_s$, $1\le s\le l$, and $k_1<k_2<\cdots <k_l$.
\end{cor}
The formula (\ref{7.4}) was stated at first in [FK], Proposition~4.2.

\begin{problem} It is an interesting task to find the Koszul dual of the 
algebra $\A_n^t$.
\end{problem}

Let us postpone a proof of Theorem~\ref{t7.1} to the end of this Section, 
and start with a review of some results obtained by V.I.~Arnold [A], 
I.M.~Gelfand and A.N.~Var\-chen\-ko [GV].

It seems to be an interesting problem to understand the connections 
between the algebra 
$\E_n^0$ and the Orlik--Solomon algebra [OS] that corresponds to the 
Coxeter hyperplane arrangement (of type $A_{n-1}$). Let us remind that 
V.I.~Arnold [A] described this algebra as the quotient of the exterior 
algebra of the vector space spanned by the $[ij]$ by the ideal generated 
by the relations (\ref{7.2}). The "even analog" of the Orlik--Solomon 
algebra of a hyperplane arrangement was introduced and studied by 
I.M.~Gelfand and A.N.~Varchenko in [GV]. For 
the Coxeter hyperplane arrangement of type $A_{n-1}$, the results of 
Theorems~5 and 6 of [GV] can be formulated as follows. 

Let ${\cal P}_n$ be 
the ring of integer--valued functions that are defined on the complement 
$M_C$ of the union of a finite set $C$ of hyperplaines in real, 
$n$--dimensional affine space, that have constant values on each 
connected component. The ring ${\cal P}_n$ has a distinguish set of 
generators, namely, the Heaviside functions of hyperplaines: for each 
hyperplane $f(x)=0$, the Heaviside function of the hyperplane $f$ takes 
the value 1 from one side (where $f(x)>0$), and 0 from the other (where 
$f(x)<0$). Any element of the ring ${\cal P}_n$ is a polynomial in the 
Heaviside functions. Now let us assume that
$$\M_C=\R^n\setminus\bigcup_{1\le i<j\le n}\{ z_i-z_j=0\} .
$$
Let $X_{ij}$ be the Heaviside function of the hyperplane $H_{ij}=\{ 
z_i-z_j=0\}\subset\R^n$, then ([GV], Theorems~5, 6, and 8)

$\bullet$ The Heaviside functions $X_{ij}$, $1\le i<j\le n$, satisfy the 
following relations:
\begin{eqnarray}
& & X_{ij}^2-X_{ij}=0, \label{7.5}\\
& & X_{ij}(X_{ik}-1)X_{jk}-(X_{jk}-1)X_{ik}(X_{ij}-1)=0, \label{7.6}
\end{eqnarray}
for any $1\le i<j<k\le n$.

$\bullet$ Let ${\cal T}_n$ be an ideal in the ring of polynomials $\Z 
[X_{ij}, \ 1\le i<j\le n]$ generated by the left--hand sides of the 
relations (\ref{7.5}) and (\ref{7.6}). Then the natural homomorphism $\Z 
[X_{ij}]/{\cal T}_n\to{\cal P}_n$ is an isomorphism.

$\bullet$ Let us define the degree of a function $p\in{\cal P}_n$ as the 
minimum of the degrees of polynomials in $\{ X_{ij}\}$ that represent 
$p$, and denote by ${\cal P}_n^k$ the subspace of functions degree not 
higher then $k$. Then there exists a noncanonical linear map
$$\pi_k~:~{\cal P}_n^k\to H^k(M_C,\Z ),
$$
which induces an isomorphism (${\cal P}_n^{-1}=\phi$)
$${\cal P}_n^k/{\cal P}_n^{k-1}~{\buildrel\sim\over\rightarrow}~H^k(M_C,\Z ), 
\ 0\le k\le n.
$$

As a corollary, we have
\begin{equation}
H({\cal P}_n;t):=\sum_{k=0}^nt^k\dim{\cal P}_n^k/{\cal P}_n^{k-1}=
\prod_{j=1}^{n-1}(1+jt). \label{7.7}
\end{equation}

Now we are ready to prove Theorem~\ref{t7.1}.

{\it Proof of Theorem~\ref{t7.1}.} Let us rewrite the relations 
(\ref{7.5}) and (\ref{7.6}) in slightly different form
\begin{eqnarray}
&& X_{ij}^2=\beta X_{ij}, \ \ 1\le i<j\le n; \label{7.8}\\
&& X_{ij}X_{jk}=X_{jk}X_{ik}+X_{ik}X_{ij}-\beta X_{ik}, 
\ \ 1\le i<j<k\le n, \label{7.9}
\end{eqnarray}
where $\beta$ is an additional deformation parameter.

We define a deformation ${\cal P}_{n,\beta}$ of the algebra ${\cal P}_n$ 
as the quotient algebra
$$\Z [X_{ij}]/{\cal T}_{n,\beta}={\cal P}_{n,\beta},
$$
where ${\cal T}_{n,\beta}$ is an ideal in the commutative ring  of 
polynomials $\Z [X_{ij}]$ generated by the relations (\ref{7.8}) and 
(\ref{7.9}).

It is clear that ${\cal P}_{n,\beta =0}\cong\A_n^0$, and ${\cal 
P}_{n,\beta}$ is a flat deformation of $\A_n^0$. Thus, we have 
%(proof of Theorem~\ref{t7.1}):
$$\prod_{j=1}^{n-1}(1+jt)=H({\cal P}_n;t)=H({\cal P}_{n,\beta};t)=
H({\cal P}_{n,\beta =0};t)=H(\A_n^0;t).
$$
\qed

%\section{Quadratic algebra $\L_{n,\beta}$ }
%\label{sec: last}
%\neweq
\vskip 0.4cm
Let us introduce a noncommutative version of the algebra 
${\cal P}_{n,\beta}$.
\begin{de}\label{d7.2} Define the algebra ${\cal L}_{n,\beta}$ (of type 
$A_{n-1}$) as the quadratic algebra over $\Z [\beta ]$ with generators $[ij]$, 
$1\le i<j\le n$, which satisfy the following relations
\begin{equation} 
\label{7.10}
\begin{array}{rl}
&[ij][jk]=[jk][ik]+[ik][ij]+\beta [ik]\,,\\[.1in]
&[jk][ij]=[ik][jk]+[ij][ik]+\beta [ik],
                        \textrm{\ \ for } i<j<k\,;\quad\ \  \\[.1in]
\end{array}
\end{equation}
\begin{equation} 
\label{7.11}
\begin{array}{rl}
&~[ij][kl]=[kl][ij],\\[.1in]
&~\textrm{whenever } \{i,j\}\cap\{k,l\}=\phi\,,
   i<j, \textrm{ and }k<l \,. ~~~~~~~~~~
\end{array}
\end{equation}
\end{de}
The algebra ${\cal L}_{n,\beta}$ is a smooth deformation of the algebra 
$\G_n$, and has the same Hilbert series:
$$H({\cal L}_{n,\beta};t):=\sum_{k=0}^{\infty}t^k\dim{\cal 
L}_{n,\beta}^{(k)}/{\cal L}_{n,\beta}^{(k-1)}=H(\G_n;t).
$$

We can define also the quotients $\L_{n,\beta}^t$, $\L_{n,\beta}^0$, and 
commutative subalgebras $\H_{n,\beta}$, $\H_{n,\beta}^t$ and 
$\H_{n,\beta}^0$ generated by the Dunkl elements $\theta_j$ (see 
(\ref{2.1})), $1\le j\le n$. 

It seems interesting to describe the 
relations in the commutative subalgebra 
$\H_{n,\beta}^t\subset\L_{n,\beta}^t$. The work is in progress and we 
hope to present our results in the nearest future.

\end{document}